\documentclass[a4paper,11pt]{article}
\usepackage{amsmath}

\newtheorem{theorem}{Theorem}[section]
\newtheorem{lemma}{Lemma}[section]

\usepackage{amssymb}
\usepackage{graphicx}
\usepackage{epstopdf}

\numberwithin{equation}{section}
\author{Alessandro De Gregorio\footnote{Dipartimento di Scienze Economiche, Aziendali e Statistiche,
Via Conservatorio 7, 20122, Milano, Italy. {\bf email}: alessandro.degregorio@unimi.it}}
\title{ Parametric estimation for planar random flights observed at discrete
times}
\begin{document}
\maketitle
\begin{abstract}
We deal with a planar random flight $\{(X(t),Y(t)),0<t\leq T\}$
observed at $n+1$ equidistant times $t_i=i\Delta_n,i=0,1,...,n$.
The aim of this paper is to estimate the unknown value of the
parameter $\lambda$, the underlying rate of the Poisson process.
The planar random flights are not markovian, then we use an
alternative argument to derive a pseudo-maximum likelihood
estimator $\hat{\lambda}$ of the parameter $\lambda$. We consider two different types of
asymptotic schemes and show the consistency, the asymptotic
normality and efficiency of the estimator proposed.  A Monte Carlo
analysis for small sample size $n$ permits us to analyze the
empirical
performance of $\hat{\lambda}$. 

A different approach permits us to introduce an alternative estimator of $\lambda$ which is consistent, asymptotically normal and asymptotically efficient without the request of other assumptions. 
\end{abstract}
\section{Introduction}
 Diffusion processes play a central role in the theory of
stochastic processes. However these models do not give a realistic
description of the real movements because the velocity is infinite
and their sample paths are nowhere differentiable. For this reason
in literature have been proposed alternative processes with finite
velocity. The first model of this type, introduced by Goldstein
(1951) and Kac (1974), is the telegraph process which describes
the motion of a particle on the real line.

The planar random flights are a possible extension in
$\mathbb{R}^2$ of the telegraph process. We consider the motion in
the plane of a particle starting at arbitrary point $(x_0,y_0)$,
moving with constant velocity $c$ and taking directions uniformly
distributed in $(0,2\pi]$. The changes of direction are governed
by a homogenous Poisson process with parameter $\lambda>0$. Let
$N(t)$ be the number of Poisson events in the interval $[0,t]$,
the position at time $t>0$ of a particle performing a random
flight is
\begin{eqnarray}
&&X(t)=x_0+c\sum_{j=1}^{N(t)+1}(s_j-s_{j-1})\cos\theta_j,\notag\\
&&Y(t)=y_0+c\sum_{j=1}^{N(t)+1}(s_j-s_{j-1})\sin\theta_j,
\end{eqnarray}
where $\theta_j $ are independent random variables uniformly
distributed in $(0,2\pi]$, while $s_j,j=1,...,n$ are the instants
at which Poisson events occur ($s_0=0$ and $s_{N(t)+1}=t).$

The distribution $p(x,y;t)$ of $(X(t),Y(t))$ is concentrated on
the disc
$$S_{ct}^2=\{(x,y):(x-x_0)^2+(y-y_0)^2\leq  c^2t^2\}.$$

If the initial direction is maintained until time $t$, the
probability density $p(x,y;t)$ possesses a singular component,
otherwise the distribution lies inside $S_{ct}^2$.

Random flights in $\mathbb{R}^2$ have been studied in Stadje
(1987), Masoliver et al. (1993), Kolesnik and Orsingher (2005). De
Gregorio and Orsingher (2007) analyze random flights in
$\mathbb{R}^d,d \geq 2$, and derive their explicit distribution in
the four-dimensional space.

The only references about the statistical inference of random
motion at finite velocity consider estimation problem for the
telegraph process. Yao (1985) estimates the state of the telegraph
process under white noise perturbation and studies performance of
nonlinear filters. Iacus (2001) considers the estimation of the
parameter $\theta$ of a non-constant rate $\lambda_\theta(t)$.
More recently, De Gregorio and Iacus (2006) introduce a
pseudo-maximum likelihood estimator and a least squares estimator
for the parameter $\lambda$ when the sample paths of the telegraph
process are observed only at equidistant discrete times. The
authors also analyze the same statistical problem for a geometric
telegraph process particularly interesting in view of financial
applications. For a telegraph process observed at discrete times,
Iacus and Yoshida (2006) study the asymptotic (i.e. the mesh
decreases to zero and the horizon interval tends to infinity)
properties of two moment estimators and propose also an estimator
consistent, asymptotically normal and efficient.

The aim of this paper is the estimation of the parameter $\lambda$
when the process $\{(X(t),Y(t)),0<t\leq T\}$ is observed at $n+1$
equidistant times $0=t_0<t_1<..<t_n=T$, where
$t_i=i\Delta_n=i\Delta,i=0,1,...,n.$ We consider two types of
asymptotic framework:

 1) $\Delta_n\rightarrow 0$ and
$n\Delta_n=T\rightarrow \infty$ as $n\rightarrow \infty$;

 2)
$n\rightarrow \infty$ with $\Delta_n$ fixed.

The statistical problem is interesting because the planar random
flights seem to be useful for ecology and biology applications. In
fact, Holmes (1993) and Holmes et al. (1994) consider these models
to represent the displacements of the animals and microorganisms
on a surface.

We note that when the planar random flight is observed
continuously, then $N(T)/T$ is the optimal estimator of the
parameter $\lambda$ and our statistical problem is equivalent to
the estimate of a whole Poisson process on $[0,T]$ (see Kutoyants
(1998)).

The process $\{(X(t),Y(t)),0<t\leq T\}$ is not markovian. Hence,
it is not possible to explicit the likelihood function of the
points observed as product of the transition densities. This fact
implies that we can not use the tools developed for the diffusion
processes (see Sorensen (1997) and Soresen (2004) for an account
of these estimation methods).

The main idea of this paper is to consider the points
$$(X(i\Delta_n)-X((i-1)\Delta_n),Y(i\Delta_n)-Y((i-1)\Delta_n))$$ as
$n$ independent copies of a random flight up to time $\Delta_n$
(which is untrue). In this manner we can build an estimating
function from which it is possible to derive a pseudo-maximum
likelihood estimator.

The paper is organized as follows. In Section 2 we describe the
random motion considered here and present some results concerning
the process $\sqrt{X^2(t)+Y^2(t)}$ (for example the moments). In
Section 3 we introduce a pseudo-likelihood function $L_n(\lambda)$
and propose the following estimator
\begin{equation}
\hat{\lambda}_n=\arg \max_{\lambda>0} L_n(\lambda).
\end{equation}
Under the asymptotic scheme 1) the estimator $\hat{\lambda}_n$ is
asymptotically normal and efficient. Alternatively under the same asymptotic scheme , we present an estimator 
asymptotically efficient without additional hypotheses. By considering
the second asymptotic framework in Section 4, we can study the convergence of the estimator $\hat{\lambda}_n$ by means of the
pseudo-likelihood ratio. In the last section, we analyze the
empirical performance of $\hat{\lambda}_n$ by means of a Monte
Carlo analysis.

\section{Planar random flights: description and some results}
We consider a particle starting at the arbitrary point $(x_0,y_0)$
of the plane $\mathbb{R}^2$, moving with constant finite speed
$c$. The initial direction is a random variable $\theta$ uniformly
distributed in $(0,2\pi]$. The changes of direction are governed
by a homogeneous Poisson process with parameter $\lambda>0$.
Therefore, when a Poisson time occurs the particle takes a new
direction uniformly distributed in $(0,2\pi]$, independently from
the previous one.

We indicate the position of the particle at time $t>0$ with the
stochastic process $(X(t),Y(t))$, which is called \emph{random
flight} in the plane. At time $t$ the particle is located in the
disc
\begin{equation}\label{disc}
 S_{ct}^2=\{(x,y):(x-x_0)^2+(y-y_0)^2\leq
c^2t^2\},
\end{equation}
with probability 1. If no Poisson event occurs the particle
reaches the circle $\partial S_{ct}^2=\{(x,y):(x-x_0)^2+(y-y_0)^2=
c^2t^2\},$ with probability
$$\mathrm{P}\{(X(t),Y(t))\in \partial S_{ct}^2\}=e^{-\lambda t}.$$

The remaining part of the distribution lies in the interior of
\eqref{disc} and represents the absolute continuous component of
the distribution
\begin{equation}
\mathrm{P}\left\{X(t)\in dx,Y(t)\in dy\right\}.
\end{equation}

We note that the random flights have trajectories which assume the
form of broken lines where the single steps have random length and
are uniformly oriented in $(0,2\pi]$. However, the total length
for any sample paths at time $t$ is $ct$.

The density law of $(X(t),Y(t))$ (see Kolesnik and Orsingher
(2005)) is equal to
\begin{eqnarray}\label{density}
p(x,y;t)&=&\frac{\lambda}{2\pi c}\frac{e^{-\lambda
t+\frac{\lambda}{c}\sqrt{c^2t^2-(x-x_0)^2-(y-y_0)^2}}}{\sqrt{c^2t^2-(x-x_0)^2-(y-y_0)^2}}\mathbf{1}_{\{(x-x_0)^2+(y-y_0)^2<c^2t^2\}}\notag\\
&&+\frac{e^{-\lambda t}}{2\pi
c}\delta(c^2t^2-(x-x_0)^2-(y-y_0)^2),
\end{eqnarray}
with $(x,y)\in S_{ct}^2$ and $\delta(\cdot),\mathbf{1}(\cdot)$
representing respectively the Dirac's delta function and the
indicator function.

Now, we present some results concerning the following process
\begin{equation}\label{distance}
R(t)=\sqrt{X^2(t)+Y^2(t)},
\end{equation}
i.e. the euclidean distance from the origin of the space
$\mathbb{R}^2$ of the position reaches by the moving particle at
time $t>0.$

The following theorem contains our first result.

\begin{theorem}
The absolute continuous component of the process $R(t),t>0$, when
$R(0)=\sqrt{x_0^2+y_0^2}$, is equal to
\begin{equation}\label{lawdistance}
f_R(r,t)=\frac{\lambda}{2\pi c}r e^{-\lambda
t}\int_0^{2\pi}\frac{e^{\frac{\lambda}{c}\sqrt{c^2t^2-r^2-x_0^2-y_0^2+2r
(x_0\cos \theta+y_0\sin\theta)}}}{\sqrt{c^2t^2-r^2-x_0^2-y_0^2+2r
(x_0\cos \theta+y_0\sin\theta)}}d\theta,
\end{equation}
with $0<r<ct$. Moreover, under the Kac condition (i.e. $c,\lambda
\rightarrow \infty$ in such a way that $\frac{c^2}{\lambda}
\rightarrow 1$), we have that \eqref{lawdistance} tends to the law
of a standard Bessel process.
\end{theorem}

\noindent{\bf Proof.} We start the proof observing that
\begin{eqnarray}\label{probdistan}
&&\mathrm{P}\left\{R(t)\leq r | R(0)=\sqrt{x_0^2+y_0^2}\right\}\\
&&=\frac{\lambda}{2\pi c}\iint_{\{(x,y):x^2+ y^2 \leq
r^2\}}\frac{e^{-\lambda
t+\frac{\lambda}{c}\sqrt{c^2t^2-(x-x_0)^2-(y-y_0)^2}}}{\sqrt{c^2t^2-(x-x_0)^2-(y-y_0)^2}}dxdy\notag\\
\notag\\&&=\frac{\lambda}{2\pi c} \int_0^r\rho
d\rho\int_0^{2\pi}\frac{e^{-\lambda
t+\frac{\lambda}{c}\sqrt{c^2t^2-(x_0-\rho\cos\theta)^2-(y_0-\rho\cos\theta)^2}}}{\sqrt{c^2t^2-(x_0-\rho\cos\theta)^2-(y_0-\rho\cos\theta)^2}}d\theta\notag\\
\notag\\&&=\frac{\lambda}{2\pi c}\int_0^r\rho
d\rho\int_0^{2\pi}\frac{e^{-\lambda
t+\frac{\lambda}{c}\sqrt{c^2t^2-\rho^2-x_0^2-y_0^2+2\rho (x_0\cos
\theta+y_0\sin\theta)}}}{\sqrt{c^2t^2-\rho^2-x_0^2-y_0^2+2\rho
(x_0\cos \theta+y_0\sin\theta)}}d\theta.\notag
\end{eqnarray}

By differentiating the probability \eqref{probdistan} with respect
to $r$, the density \eqref{lawdistance} emerges.

In order to prove the second part of the theorem, we rewrite the
density \eqref{lawdistance} in the following form
\begin{eqnarray}\label{exbessel}
f_R(r,t)&=&\frac{\lambda}{2\pi c}r \int_0^{2\pi}\frac{e^{-\lambda
t+\lambda t\sqrt{1-\frac{r^2+x_0^2+y_0^2-2r (x_0\cos
\theta+y_0\sin\theta)}{c^2t^2}}}}{ct\sqrt{1-\frac{r^2+x_0^2+y_0^2-2r
(x_0\cos \theta+y_0\sin\theta)}{c^2t^2}}}d\theta\notag\\
&=&\frac{\lambda}{2\pi c^2 t}r
\int_0^{2\pi}\frac{e^{-\frac{r^2+x_0^2+y_0^2-2r (x_0\cos
\theta+y_0\sin\theta)}{2\frac{c^2}{\lambda}t}-...}}{\sqrt{1-\frac{r^2+x_0^2+y_0^2-2r
(x_0\cos \theta+y_0\sin\theta)}{c^2t^2}}}d\theta.
\end{eqnarray}
In the last step we have used the expansion
$\sqrt{1-w}=1-\frac{w}{2}-\frac{w^2}{8}-...,$ which is absolutely
convergent for $|w|<1$.

From \eqref{exbessel}, under the Kac condition, we obtain the
following limit
\begin{eqnarray}
\lim_{\substack{\lambda,c \rightarrow \infty\\
c^2/\lambda \rightarrow 1}}f_R(r,t)&=&\frac{r}{2\pi
t}\int_0^{2\pi}e^{-\frac{r^2}{2t}}e^{-\frac{x_0^2+y_0^2}{2t}}e^{\frac{r
(x_0\cos \theta+y_0\sin\theta)}{t}}d\theta\notag\\
&=&
\frac{r}{t}e^{-\frac{r^2}{2t}}e^{-\frac{x_0^2+y_0^2}{2t}}I_0\left(\frac{r\sqrt{x_0^2+y_0^2}}{t}\right),\label{bessel}
\end{eqnarray}
by means of the well-known integral representation of the Bessel
function
$$I_0(x\sqrt{\alpha^2+\beta^2})=\frac{1}{2\pi}\int_0^{2\pi}e^{x(\alpha\cos\theta+\beta\sin\theta)}d\theta.$$

Expression \eqref{bessel} represents the density function of
Bessel process $\sqrt{B_1^2(t)+B_2^2(t)}$, where $B_1$ and $B_2$
are two
independent standard Brownian motion.$\hfill\square$\\

From Theorem 2.1 we note that for $(x_0,y_0)=(0,0)$, the complete
distribution of $R(t)$ becomes
\begin{eqnarray}\label{besselorigin}
p_R(r,t)&=&\frac{\lambda}{c}\frac{r\exp\{-\lambda
t+\frac{\lambda}{c}\sqrt{c^2t^2-r^2}\}}{\sqrt{c^2t^2-r^2}}\mathbf{1}_{\{0<r<ct\}}\\
&&+\frac{re^{-\lambda t}}{ct}\delta(c^2t^2-r^2).\notag
\end{eqnarray}

We note that the density \eqref{besselorigin} coincides with
formula (7) in Kolesnik and Orsingher (2005), when we ignore the
angular component.

By taking into account the probability law \eqref{besselorigin},
we are able to derive the moments of $R(t)$.

\begin{theorem}
Let $(x_0,y_0)=(0,0)$ and $p\geq 1$, we have that
\begin{equation}\label{moments}
\mathrm{E}R^p(t)=(ct)^pe^{-\lambda
t}\left\{\sqrt{\pi}\left(\frac{2}{\lambda
t}\right)^\frac{p-1}{2}\Gamma\left(\frac{p+1}{2}\right)I_\frac{p+1}{2}(\lambda
t)+1\right\}.
\end{equation}
\end{theorem}
\noindent{\bf Proof.} In view of \eqref{besselorigin}, we can
write
\begin{eqnarray}\label{calmom}
\mathrm{E}R^p(t)=\frac{\lambda}{c}e^{-\lambda t}\int_0^{ct}r^{p+1}
\frac{e^{\frac{\lambda}{c}\sqrt{c^2t^2-r^2}}}{\sqrt{c^2t^2-r^2}}dr+(ct)^pe^{-\lambda
t}.
\end{eqnarray}

Now, we work out the integral in \eqref{calmom}. Hence
\begin{eqnarray}
&&\int_0^{ct}r^{p+1}
\frac{e^{\frac{\lambda}{c}\sqrt{c^2t^2-r^2}}}{\sqrt{c^2t^2-r^2}}dr\notag\\
&&=\sum_{k=0}^\infty
\frac{1}{k!}\left(\frac{\lambda}{c}\right)^k\int_0^{ct}r^{p+1}(c^2t^2-r^2)^{\frac{k-1}{2}}dr=(r=ct\sqrt{y})\notag\\
&&=\sum_{k=0}^\infty
\frac{1}{k!}\left(\frac{\lambda}{c}\right)^k\frac{(ct)^{p+k+1}}{2}\int_0^1y^\frac{p}{2}(1-y)^\frac{k-1}{2}dy\notag\\
&&=\sum_{k=0}^\infty
\frac{1}{k!}\left(\frac{\lambda}{c}\right)^k\frac{(ct)^{p+k+1}}{2}
\frac{\Gamma\left(\frac{p}{2}+1\right)\Gamma\left(\frac{k+1}{2}\right)}{\Gamma\left(\frac{k+1}{2}+\frac{p}{2}+1\right)}\notag\\
&&=\sqrt{\pi}\Gamma\left(\frac{p}{2}+1\right)\sum_{k=0}^\infty
\frac{1}{k!}\left(\frac{\lambda
t}{2}\right)^k\frac{(ct)^{p+1}\Gamma(k)}{\Gamma\left(\frac{k+1}{2}+\frac{p}{2}+1\right)\Gamma\left(\frac{k}{2}\right)}=(k=2m)\notag\\
&&=\sqrt{\pi}\Gamma\left(\frac{p}{2}+1\right)(ct)^{p+1}\sum_{m=0}^\infty\left(\frac{\lambda
t}{2}\right)^{2m}
\frac{1}{2m\Gamma(m)\Gamma\left(m+\frac{p+1}{2}+1\right)}\notag\\
&&=\frac{\sqrt{\pi}}{2}\Gamma\left(\frac{p}{2}+1\right)(ct)^{p+1}\left(\frac{2}{\lambda
t}\right)^\frac{p+1}{2}\sum_{m=0}^\infty\frac{1}{m!}\left(\frac{\lambda
t}{2}\right)^{2m+\frac{p+1}{2}}\frac{1}{\Gamma\left(m+\frac{p+1}{2}+1\right)}\notag\\
&&=\frac{\sqrt{\pi}}{2}\Gamma\left(\frac{p}{2}+1\right)(ct)^{p+1}\left(\frac{2}{\lambda
t}\right)^\frac{p+1}{2}I_\frac{p+1}{2}(\lambda
t).\label{calmom.bis}
\end{eqnarray}

By inserting \eqref{calmom.bis} into \eqref{calmom} we obtain the
result \eqref{moments}.$\hfill\square$\\

{\bf Remark 2.1} We observe that
\begin{equation}
\lim_{\lambda\rightarrow \infty}\mathrm{E}R^p(t)=0.
\end{equation}

In other words, if $\lambda$ grows to infinity the changes of
direction increase and
consequently the distance from the origin decreases.\\

 {\bf Remark 2.2} We derive from \eqref{moments}
the mean value  of $R(t)$
\begin{equation}
\mathrm{E}R(t)=ct e^{-\lambda t}\{\sqrt{\pi}I_1(\lambda t)+1\}.
\end{equation}

In the particular case $p=2$, we can write the square mean in
terms of simple function. In fact, by means of the following
relationship
$$I_\frac{3}{2}(x)=\sqrt{\frac{2}{\pi x^3}}(x\cosh x-\sinh x),$$
we get that
\begin{eqnarray}
\mathrm{E}R^2(t)&=&(ct)^2e^{-\lambda
t}\left\{\sqrt{\pi}\frac{\lambda t\cosh \lambda t-\sinh\lambda
t}{(\lambda t)^2}+1\right\}\notag\\
&=&(ct)^2\left\{\sqrt{\pi}\frac{\lambda t(1+e^{-2\lambda
t})-1+e^{-2\lambda t}}{2(\lambda t)^2}+e^{-\lambda t}\right\}.
\end{eqnarray}

\section{Parametric estimation for planar random flights}

We assume that the planar random flight $\{(X(t),Y(t)),0<t\leq
T\}$, with $(X(0),Y(0))=(0,0)$, is observed only at $n+1$
equidistant discrete times $0=t_0<t_1<...<t_n=T,$ where
$t_i=i\Delta_n=i\Delta,\,i=0,1,...,n.$ We use the following
notation to simplify the formulas:
$(X(t_i),Y(t_i))=(X(i\Delta_n),Y(i\Delta_n))=(X_i,Y_i).$

The interest is the estimation of the parameter $\lambda$ whilst
the velocity $c$ is assumed to be known. In other words, we want
to estimate the rate of change of a microorganism which performs a
planar random flight, when we are able to observe its position
only at discrete times.

The estimation of $c$ is an uninteresting problem. In fact, if in
the interval $(i\Delta_{n-1},i\Delta_{n}]$ there are not changes
of direction, then
$(X_i-X_{i-1})^2+(Y_i-Y_{i-1})^2=c^2\Delta_n^2$. If $\Delta_n$ is
suitable small, there is high probability of observing
$N(t_i)-N(t_{i-1})=0$ and $c$ can be calculated without error.

Analogously to the telegraph process, the random flights are not
markovian. For this reason we cannot write the explicit likelihood
of the process in the form of product of transition densities as
well as  for diffusion processes. Therefore, we need an
alternative argument in the spirit of the paper by De Gregorio and
Iacus (2006).

We define a pseudo-likelihood function as follows. By taking into
account the distribution \eqref{density}, we introduce the
following data dependent function
\begin{eqnarray}\label{applik}
L_n(\lambda)
&=&L_n(\lambda|(X_0,Y_0),(X_1,Y_1),...,(X_n,Y_n))\\
&=&\prod_{i=1}^n p((X_i,Y_i),\Delta_n; (X_{i-1},Y_{i-1}),t_{i-1})\notag\\
 &=&\prod_{i=1}^n\left\{\frac{\lambda}{2\pi
c}\frac{\exp\{-\lambda
\Delta_n+\frac{\lambda}{c}\sqrt{u_{n,i}}\}}{\sqrt{u_{n,i}}}
\mathbf{1}_{\{u_{n,i}>0\}}+\frac{e^{-\lambda \Delta_n}}{2\pi
c}\delta(u_{n,i}=0)\right\},\notag
\end{eqnarray}
where
$u_{n,i}=u_n((X_i,Y_i),(X_{i-1},Y_{i-1}))=c^2\Delta_n^2-(X_{i}-X_{i-1})^2-(Y_{i}-Y_{i-1})^2.$

The transition densities $p((X_i,Y_i),\Delta_n;
(X_{i-1},Y_{i-1}),t_{i-1})$ appearing in \eqref{applik} represent
the distribution of a random flight in $\mathbb{R}^2$, initially
located at the point $(X_{i-1},Y_{i-1})$ at time $t_{i-1}$, which
reaches the position $(X_{i},Y_{i})$ at the instant $t_i$. The
function \eqref{applik} is indeed the joint law of the points
$(X_i-X_{i-1},Y_i-Y_{i-1})$, which are considered as if they were
$n$ independent copies of the process $(X(\Delta_n),Y(\Delta_n))$
(i.e. the process $(X(t),Y(t))$ up to time $\Delta_n$).

\begin{figure}
\begin{center}
\includegraphics[angle=0,width=0.8\textwidth]{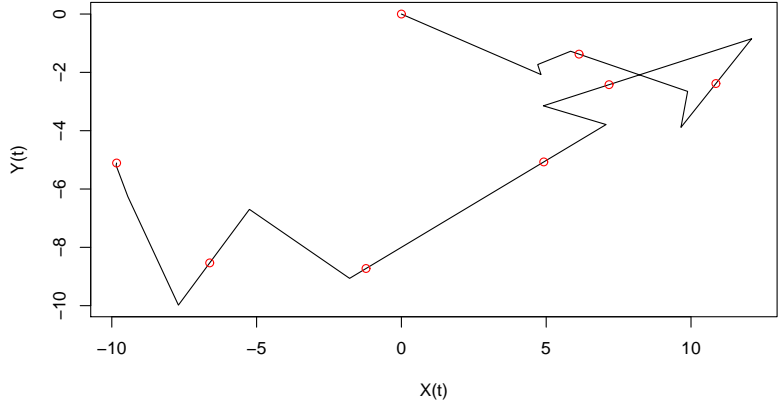}
\caption{ Discrete time sampling of the planar random flight. For
this sample path $n=7$ and $n^+=6$.} \label{fig}
\end{center}
\end{figure}

The pseudo-likelihood function \eqref{applik} is equivalent to
\begin{eqnarray}\label{applik.bis}
L_n(\lambda) &=&\left(\frac{e^{-\lambda \Delta_n}}{2\pi
c}\right)^{(n-n^+)}\prod_{i=1}^{n^+}\left\{\frac{\lambda}{2\pi
c}\frac{\exp\{-\lambda
\Delta_n+\frac{\lambda}{c}\sqrt{u_{n,i}}\}}{\sqrt{u_{n,i}}} \right\}\notag\\
&=&\frac{e^{-\lambda n\Delta_n }}{(2\pi
c)^{n}}\frac{\lambda^{n^+}\exp\left\{\frac{\lambda}{c}\sum_{i=1}^{n^+}\sqrt{u_{n,i}}\right\}}{\prod_{i=1}^{n^+}\sqrt{u_{n,i}}},
\end{eqnarray}
where $n^+$ is the number of the planar random flights with at
least one change of direction.

{\bf Remark 3.1} In the expression \eqref{applik.bis}, the factor
$\left(\frac{e^{-\lambda \Delta_n}}{2}\right)^{n-n^+}$ concerns
the singular part of the densities $p((X_i,Y_i),\Delta_n;
(X_{i-1},Y_{i-1}),t_{i-1})$, while the product represents the
absolutely continuous components of the distributions of the
random flights. Note that for increasing values of $\lambda$, the
absolutely continuous component of \eqref{applik.bis} has a bigger
weight than the singular component; viceversa for small values of
$\lambda$. Figure \ref{fig} shows how the two components of the
function $L_n(\lambda)$ emerge
for this scheme of observation.\\

To define our estimator, we borrow the approach based on
estimating functions. Formula \eqref{applik.bis} yields
\begin{eqnarray}\label{score}
F_n(\lambda)=\frac{d}{d\lambda}\log L_n(\lambda) =-n\Delta_n
+\frac{1}{c}\sum_{i=1}^{n^+}\sqrt{u_{n,i}} +\frac{n^+}{\lambda},
\end{eqnarray}
and the estimator is obtained by solving $F_n(\lambda)=0$.

By taking into account the function \eqref{score} we derive the
following pseudo-maximum likelihood estimator for the parameter
$\lambda$
\begin{eqnarray}\label{maxlika}
\hat{\lambda}_n&=&\arg \max_{\lambda>0}
L_n(\lambda)\notag\\
&=&\frac{cn^+}{cn\Delta_n-\sum_{i=1}^{n^+} \sqrt{u_{n,i}}}.
\end{eqnarray}

It is easy to see that $\frac{d^2}{d\lambda^2}\log L_n(\lambda)<0$
and the uniqueness of the estimator \eqref{maxlika} holds. However
$\hat{\lambda}_n$ is not a true maximum likelihood estimator
because $L_n(\lambda)$ is not a true likelihood function. In the
section 4 we will analyze the empirical performance of the
estimator \eqref{maxlika} for small sample size.

By considering the following asymptotic framework
$\Delta_n\rightarrow 0$ and $n\Delta_n=T\rightarrow \infty$ as
$n\rightarrow\infty$, we provide the next result for the estimator
\eqref{maxlika}.
\begin{theorem}
Let  $\Delta_n\rightarrow 0$ and $n\Delta_n=T\rightarrow \infty$
 as $n \rightarrow \infty$, then $\hat{\lambda}_n$ is consistent, asymptotically normal and
 efficient.
 \end{theorem}
{\bf Proof.} If $\Delta_n \rightarrow 0$ and $n\Delta_n=T$ as
$n\rightarrow \infty$, we have that $u_{n,i}\rightarrow 0$ and
$n^+\rightarrow N(T)$, where $N(T)$ represents the number of
changes of direction occurred in the interval $[0,T]$ when the
whole trajectory is observed. Thus, $\hat{\lambda}_n$ tends to the
maximum likelihood estimator of a homogeneous Poisson process
\begin{eqnarray}
\frac{N(T)}{T}.
\end{eqnarray}

Therefore, under the condition $n\Delta_n=T\rightarrow \infty$,
the pseudo-maximum likelihood estimator is consistent,
asymptotically
normal and efficient (see Kutoyants (1998)).$\hfill\square$\\

We introduce an alternative estimator for the parameter $\lambda$ by means of the distances
$$\eta_i=\sqrt{(X_i-X_{i-1})^2+(Y_i-Y_{i-1})^2}.$$

By setting 
\begin{equation*}
G_n=\frac{1}{n\Delta_n}\sum_{i=1}^n\mathbf{1}_{\{\eta_i<c\Delta_n\}}
=\frac{1}{n\Delta_n}\sum_{i=1}^n\mathbf{1}_{\{N([t_{i-1},t_i))\geq
1\}}.
\end{equation*}
we define the following unbiased estimator
\begin{equation}\label{alest}
\dot{\lambda}_n=-\frac{1}{\Delta_n}\log\left(1-\Delta_nG_n\right).
\end{equation}

The advantage of $\dot{\lambda}_n$ is that we are able to derive the asymptotic properties without assumptions on the points $(X_i-X_{i-1},Y_i-Y_{i-1})$.
\begin{theorem}
For $n\Delta_n=T\rightarrow\infty$ and $\Delta_n	\rightarrow 0$ as $n\rightarrow \infty$ the estimator \eqref{alest} is consistent, asymptotically normal and efficient
\begin{equation}
\sqrt{n\Delta_n}(\dot{\lambda}_n-\lambda)\overset{d} {\rightarrow} N(0,\lambda).
\end{equation}
\end{theorem}
{\bf Proof.}  We replace the steps contained in the proof presented
by Iacus and Yoshida (2006) for the telegraph process.

 First of all we prove consistency and asymptotic normality of $G_n$. We observe that
 $$\mathrm{E}(G_n)=\frac{1-e^{-\lambda\Delta_n}}{\Delta_n}=\lambda+\frac{1}{2}\lambda^2\Delta_n+o(\Delta_n^2)\rightarrow \lambda$$
 and consistency immediately follows. Now we show the asymptotic normality of $G_n$ and consider to this scope the quantity
 
 \begin{eqnarray*}
U_n&=&\sqrt{n\Delta_n}(G_n-\mathrm{E}(G_n))\\
&=&\frac{1}{\sqrt{n\Delta_n}}\sum_{i=1}^n\left[\mathbf{1}_{\{\eta_i<c\Delta_n\}}-
\mathrm{E}(\mathbf{1}_{\{\eta_i<c\Delta_n\}})\right]\\
&=&\frac{1}{\sqrt{n\Delta_n}}\sum_{i=1}^n\left[\mathbf{1}_{\{N([t_{i-1},t_i))\geq
1\}}-
(1-e^{-\lambda\Delta_n})\right]\\
&=&\sum_{i=1}^n\alpha_i
\end{eqnarray*}
where
$$\alpha_i=\frac{1}{\sqrt{n\Delta_n}}\left[\mathbf{1}_{\{N([t_{i-1},t_i))\geq
1\}}-
(1-e^{-\lambda\Delta_n})\right].$$

It's clear that $\mathrm{E}(\alpha_i)=0$ and $\mathrm{E}(U_n)=0.$ Moreover
\begin{eqnarray*}
\mathrm{Var}(\alpha_i)&=&\frac{1}{n\Delta_n}\mathrm{Var}(\mathbf{1}_{\{N([t_{i-1},t_i))\geq
1\}})\\
&=&\frac{1}{n\Delta_n}\left\{\mathrm{E}(\mathbf{1}_{\{N([t_{i-1},t_i))\geq
1\}})-(\mathrm{E}(\mathbf{1}_{\{N([t_{i-1},t_i))\geq
1\}}))^2\right\}\\
&=& \frac{1}{n\Delta_n}\left\{e^{-\lambda\Delta_n}-e^{-2\lambda\Delta_n}\right\}\\
&=&\frac{1}{n}\left\{\lambda+o(1)\right\},
\end{eqnarray*}
therefore
\begin{equation}
\mathrm{Var}(U_n)=\lambda+o(1).
\end{equation}

The variables $\alpha_i$ are independent and the Lindeberg condition is true, i.e.
\begin{equation}\label{lind}
\sum_{i=1}^n\mathrm{E}\left\{\mathbf{1}_{\{N([t_{i-1},t_i))\geq
1\}}\alpha_i^2\right\}\rightarrow 0,
\end{equation}
because for large $n$ it holds true that $|\alpha_i| \leq  \frac{1}{\sqrt{n\Delta_n}}.$ From condition 
\eqref{lind} follows that
\begin{equation}
U_n\overset{d} {\rightarrow} N(0,\lambda).
\end{equation}

Finally, we can prove the asymptotic normality of $\dot{\lambda}_n$.  Since
$$f(w)=-\frac{1}{\Delta_n}\log(1-w\Delta_n),\quad f'(w)=\frac{1}{1-w\Delta_n},$$
and 
$$\dot{\lambda}_n=f(G_n),\quad \lambda=f\left(\mathrm{E}(
G_n)\right),$$
then, by so-called $\delta$-method, we have that
\begin{eqnarray*}
\sqrt{n\Delta_n}(\dot{\lambda}_n-\lambda)&=&\sqrt{n\Delta_n}(f(G_n)-f(\mathrm{E}(
G_n)))\\
&=&\sqrt{n\Delta_n}(G_n-\mathrm{E}(
G_n))f'(\lambda)+o_p(\sqrt{n\Delta_n}|G_n-\mathrm{E}(
G_n)|)\\
&=&\sqrt{n\Delta_n}(G_n-\mathrm{E}(
G_n))\frac{1}{1-\lambda\Delta_n}+o_p(\sqrt{n\Delta_n}|G_n-\mathrm{E}(
G_n)|),
\end{eqnarray*}
hence for $n\Delta_n=T\rightarrow\infty$ and $\Delta_n	\rightarrow 0$ as $n\rightarrow \infty$, we obtain that
\begin{equation*}
\sqrt{n\Delta_n}(\dot{\lambda}_n-\lambda)\overset{d} {\rightarrow} N(0,\lambda).
\end{equation*}
$\hfill\square$

\section{Large sample properties for the pseudo-maximum
likelihood estimator}

To analyze the properties of the estimator \eqref{maxlika} as
$n\rightarrow \infty$ with $\Delta_n$ fixed (large sample scheme),
we use the tools of asymptotic theory of statistical estimation
presented in Ibragimov and Has'minskii (1981).

Let us assume the velocity $c$ known and
$\lambda\in(\lambda_1,\lambda_2)=\Theta$ with
$0<\lambda_1<\lambda_2<\infty$.  We need to introduce another
hypothesis: the distance $\Delta_n$ between two consecutive
instants $t_i,i=0,1,...,n,$ is such that the following condition
holds
\begin{equation}\label{cond}
\mathrm{P}_\lambda\left\{\left(X(\Delta_n),Y(\Delta_n)\right)\in
\mathrm{int} \,S_{c\Delta_n}^2\right\}=1,
\end{equation}
where
$\mathrm{int}\,S_{c\Delta_n}^2=\{(x,y):x^2+y^2<c^2\Delta_n^2\}$.

In other words between the points $(X_{i-1},Y_{i-1})$ and
$(X_i,Y_i),i=1,..,n,$ the planar random flights have at least one
change of direction (or equivalently
$\mathrm{P_\lambda}\{N(i\Delta_n)=0\}=0,\,i=1,...,n$). In general
it is obvious that for increasing values of $\lambda$, the minimum
value of $\Delta_n$ satisfying the condition \eqref{cond}
decreases.

Immediately, from \eqref{cond} follows that the singular part of
\eqref{applik.bis} vanishes. In fact we have that
\begin{eqnarray}\label{cap.4:lik.con}
\widetilde{L}_n(\lambda)&=&\prod_{i=1}^n\left\{\frac{\lambda}{2\pi
c}\frac{\exp\{-\lambda
\Delta_n+\frac{\lambda}{c}\sqrt{u_{n,i}}\}}{\sqrt{u_{n,i}}}\right\}\notag\\
&=&\left(\frac{\lambda}{2\pi c}\right)^n\frac{\exp\left\{-\lambda
n
\Delta_n+\frac{\lambda}{c}\sum_{i=1}^{n}\sqrt{u_{n,i}}\right\}}{\prod_{i=1}^{n}\sqrt{u_{n,i}}},
\end{eqnarray}
while the pseudo-maximum likelihood estimator \eqref{maxlika}
becomes
\begin{eqnarray}\label{cap.4:maxlikmod}
\widetilde{\lambda}_n=\frac{cn}{cn\Delta_n-\sum_{i=1}^{n}
\sqrt{u_{n,i}}}.
\end{eqnarray}

We start our analysis observing that the Radon-Nikodym theorem
yields
\begin{equation}\label{cap.4:den}
p((X_i,Y_i),\Delta_n;
(X_{i-1},Y_{i-1}),t_{i-1})=\frac{d\mathrm{P}_{\lambda}}{d\mu}=\left(\frac{\lambda}{2\pi
c}\right)\frac{\exp\left( -\lambda
t+\frac{\lambda}{c}\sqrt{u_{n,i}}\right)}{\sqrt{u_{n,i}}},
\end{equation}
where $\mu$ is the Lebesgue measure in the plane. Thus, we can
indicate $\widetilde{L}_n(\lambda)$  as follows
\begin{equation}
\widetilde{L}_n(\lambda)=\frac{d \mathbf{P}_\lambda^n}{d\mu^n},
\end{equation}
where $\mathbf{P}_\lambda^n$ represents the joint probability distribution
of $n$ independent copies of a planar random flight up to time
$\Delta_n$.

It's appropriate to remark that we are presenting in this section
results valid only for the parametric model
$\{\mathbf{P}_\lambda^n,\lambda\in \Theta\}$, i.e. the model deriving from the assumption of i.i.d. observations.

 To
simplify the formulas we write $p((X_i,Y_i),\Delta_n;
(X_{i-1},Y_{i-1}),t_{i-1})=p(\lambda)$. Our first result is the
following theorem.
\begin{theorem}
Let $\mathcal{E}^{n}$ be the experiment generated by $n$
independent observations of $(X(\Delta_n),Y(\Delta_n))$. Then
$\mathcal{E}^n$ is regular with Fisher's information equal to
\begin{equation}\label{fisher}
I_n(\lambda)=\frac{n}{\lambda^2}.
\end{equation}
\end{theorem}
\textbf{Proof.} By considering the definition of regular
experiment presented in Ibragimov and Has'minskii (1981), page 65,
we must prove that $\sqrt{p(\lambda)}$ is differentiable in
$\mathbf{L}_2$ (the space of the square integrable functions) with
continuous derivative in $\mathbf{L}_2$
\begin{equation}
\psi(\lambda)=\frac{\sqrt{p(\lambda)}}{2}\left(-\Delta_n+\frac{1}{c}\sqrt{u_{n,i}}+\frac{1}{\lambda}\right).
\end{equation}

By setting $g(\lambda)=\sqrt{p(\lambda)}$ we get that
\begin{eqnarray}\label{diffcon}
&&\iint_{S_{c\Delta_n}^2}(g(\lambda+h)-g(\lambda)-h\psi(\lambda))^2dxdy\\
&&=\mathrm{E}_\lambda\left\{\frac{g(\lambda+h)}{g(\lambda)}-1-
\frac{h}{2}\left(-\Delta_n+\frac{1}{c}\sqrt{u_{n,i}}+\frac{1}{\lambda}\right)\right\}^2\notag\\
&&=\mathrm{E}_\lambda\left\{e^{-\frac{h\Delta_n}{2}+\frac{h}{2c}\sqrt{u_{n,i}}+\log\sqrt{\frac{\lambda+h}{\lambda}}}-1-
\frac{h}{2}\left(-\Delta_n+\frac{1}{c}\sqrt{u_{n,i}}+\frac{1}{\lambda}\right)\right\}^2.\notag
\end{eqnarray}
Now, by observing that
\begin{equation}\label{expansion}
e^{-\frac{h\Delta_n}{2}+\frac{h}{2c}\sqrt{u_{n,i}}+\log\sqrt{\frac{\lambda+h}{\lambda}}}=1+\frac{h}{2}\left(-\Delta_n+\frac{1}{c}\sqrt{u_{n,i}}+\frac{1}{\lambda}\right)+o(h),
\end{equation}
 we obtain
\begin{eqnarray*}
\iint_{S_{c\Delta_n}^2}(g(\lambda+h)-g(\lambda)-h\psi(\lambda))^2dx
dy=o(|h|^2).
\end{eqnarray*}

The continuity of $\psi(\lambda)$ is shown by means of the
dominated convergence theorem.

To complete the proof we verify that $\mathcal{E}^n$ possesses
finite Fisher's information $I_n(\lambda)$ for any $\lambda\in
\Theta$. Clearly $I_n(\lambda)=nI(\lambda)$, where $I(\lambda)$
represents Fisher's information of a single experiment. Thus, we
can write
\begin{eqnarray}\label{int}
I(\lambda)&=&4\iint_{S_{c\Delta_n}^2}|\psi(\lambda)|^2dxdy\\
&=&\frac{\lambda}{2\pi c}\iint_{S_{c\Delta_n}^2}\frac{ e^{-\lambda
\Delta_n+\frac{\lambda}{c}\sqrt{u_{n,i}}}}{\sqrt{u_{n,i}}}\left(\frac{1}{\lambda}-\Delta_n+\frac{1}{c}\sqrt{u_{n,i}}\right)^2dxdy\notag\\
&=&\frac{\lambda}{2\pi c}\iint_{S_{c\Delta_n}^2}\frac{e^{-\lambda
\Delta_n+\frac{\lambda}{c}\sqrt{u_{n,i}}}}{\sqrt{u_{n,i}}}\notag\\
&&\qquad\times\left(\left(\frac{1}{\lambda}-\Delta_n\right)^2+\frac{u_{n,i}}{c^2}+\frac{2}{c}\left(\frac{1}{\lambda}-\Delta_n\right)
\sqrt{u_{n,i}}\right)dxdy\notag\\
&=&\frac{\lambda}{2\pi
c}\int_0^{c\Delta_n}d\rho\int_0^{2\pi}d\theta\frac{\rho
e^{-\lambda
\Delta_n+\frac{\lambda}{c}\sqrt{c^2\Delta_n^2-\rho^2}}}{\sqrt{c^2\Delta_n^2-\rho^2}}\notag\\
&&\qquad\times\left(\left(\frac{1}{\lambda}-\Delta_n\right)^2+\frac{c^2\Delta_n^2-\rho^2}{c^2}+\frac{2}{c}\left(\frac{1}{\lambda}-\Delta_n\right)
\sqrt{c^2\Delta_n^2-\rho^2}\right),\notag
\end{eqnarray}
where in the last step we have used the transformation in polar
coordinates $x=x_0+\rho\cos\theta,y=y_0+\rho\sin\theta.$

To obtain the explicit value of \eqref{int} we calculate the
following three integrals
\begin{eqnarray*}
\mathcal{I}_1&=&\frac{\lambda}{2\pi
c}\left(\frac{1}{\lambda}-\Delta_n\right)^2\int_0^{c\Delta_n}d\rho\int_0^{2\pi}d\theta\frac{\rho
e^{-\lambda
\Delta_n+\frac{\lambda}{c}\sqrt{c^2\Delta_n^2-\rho^2}}}{\sqrt{c^2\Delta_n^2-\rho^2}}\\
&=&\frac{\lambda}{
c}\left(\frac{1}{\lambda}-\Delta_n\right)^2e^{-\lambda
\Delta_n}\int_0^{c\Delta_n}d\rho\frac{\rho
e^{\frac{\lambda}{c}\sqrt{c^2t^2-\rho^2}}}{\sqrt{c^2\Delta_n^2-\rho^2}}\\
&=& \left(\frac{1}{\lambda}-\Delta_n\right)^2e^{-\lambda
\Delta_n}\left(-e^{\frac{\lambda}{c}\sqrt{c^2\Delta_n^2-\rho^2}}\right)\Big|_{\rho=0}^{\rho=c\Delta_n} \\
&=& \left(\frac{1}{\lambda}-\Delta_n\right)^2\left(1-e^{-\lambda
\Delta_n}\right),
\end{eqnarray*}
\begin{eqnarray*}
\mathcal{I}_2&=&\frac{\lambda}{ 2\pi
c^3}\int_0^{c\Delta_n}d\rho\int_0^{2\pi}d\theta\rho
\sqrt{c^2\Delta_n^2-\rho^2}e^{-\lambda
\Delta_n+\frac{\lambda}{c}\sqrt{c^2\Delta_n^2-\rho^2}}\\
&=&\frac{\lambda}{ c^3}e^{-\lambda \Delta_n}\int_0^{c\Delta_n}\rho
\sqrt{c^2\Delta_n^2-\rho^2}e^{\frac{\lambda}{c}\sqrt{c^2\Delta_n^2-\rho^2}}d\rho=(z=\sqrt{c^2\Delta_n^2-\rho^2})\\
&=&\frac{\lambda}{ c^3}e^{-\lambda \Delta_n}\int_0^{c\Delta_n}z^2
e^{\frac{\lambda}{c}z}dz=e^{-\lambda \Delta_n}\left\{\frac{1}{
c^2}z^2e^{\frac{\lambda}{c}z}\Big|_{z=0}^{z=c\Delta_n}-\frac{2}{c^2}
\int_0^{c\Delta_n}ze^{\frac{\lambda}{c}z}dz\right\}\\
&=&e^{-\lambda \Delta_n}\left\{\Delta_n^2e^{\lambda
\Delta_n}-\frac{2}{c\lambda}ze^{\frac{\lambda}{c}
z}\Big|_{z=0}^{z=c\Delta_n}+\frac{2}{c\lambda}
\int_0^{c\Delta_n}e^{\frac{\lambda}{c}z}dz\right\}\\
&=&e^{-\lambda \Delta_n}\left\{\Delta_n^2e^{\lambda
\Delta_n}-\frac{2\Delta_n}{\lambda}e^{\lambda\Delta_n}+\frac{2}{\lambda^2}
e^{\frac{\lambda}{c}z}\Big|_{z=0}^{z=c\Delta_n}\right\}\\
&=&\Delta_n^2-\frac{2\Delta_n}{\lambda}+\frac{2}{\lambda^2}\left(1-e^{-\lambda
\Delta_n}\right),
\end{eqnarray*}
\begin{eqnarray*}
\mathcal{I}_3&=&\frac{\lambda}{ \pi
c^2}\left(\frac{1}{\lambda}-\Delta_n\right)\int_0^{c\Delta_n}
d\rho\int_0^{2\pi}d\theta \rho e^{-\lambda
\Delta_n+\frac{\lambda}{c}\sqrt{c^2\Delta_n^2-\rho^2}}\\
&=&\frac{2\lambda}{
c^2}\left(\frac{1}{\lambda}-\Delta_n\right)e^{-\lambda
\Delta_n}\int_0^{c\Delta_n}\rho e^{\frac{\lambda}{c}\sqrt{c^2\Delta_n^2-\rho^2}}d\rho=(z=\sqrt{c^2\Delta_n^2-\rho^2})\\
&=&\frac{2\lambda}{
c^2}\left(\frac{1}{\lambda}-\Delta_n\right)e^{-\lambda
\Delta_n}\int_0^{c\Delta_n}ze^{\frac{\lambda}{c}z}dz\\
&=&\left(\frac{1}{\lambda}-\Delta_n\right)e^{-\lambda
\Delta_n}\left\{\frac{2}{c}z e^{\frac{\lambda}{c}z}\Big|_{z=0}^{z=c\Delta_n}-\frac{2}{c}\int_0^{c\Delta_n}e^{\frac{\lambda}{c}z}dz\right\}\\
&=&\left(\frac{1}{\lambda}-\Delta_n\right)e^{-\lambda
\Delta_n}\left\{2\Delta_n e^{\lambda \Delta_n}-\frac{2}{\lambda}e^{\frac{\lambda}{c}z }\Big|_{z=0}^{z=c\Delta_n}\right\}\\
&=&2\left(\frac{1}{\lambda}-\Delta_n\right)\left(\Delta_n-\frac{1}{\lambda}\left(1-e^{-\lambda
\Delta_n}\right)\right).
\end{eqnarray*}

Putting together $\mathcal{I}_1,\mathcal{I}_2,\mathcal{I}_3$ we
have that
\begin{eqnarray*}
I(\lambda)&=&\mathcal{I}_1+\mathcal{I}_2+\mathcal{I}_3\\
&=&\left(1-e^{-\lambda
\Delta_n}\right)\left(\frac{1}{\lambda^2}+\Delta_n^2\right)-\Delta_n^2\\
&=&\frac{1}{\lambda^2}\left(1-e^{-\lambda
\Delta_n}\left(1+\lambda^2\Delta_n^2\right)\right).
\end{eqnarray*}

By assumption $\mathrm{P}_\lambda\{N(\Delta_n)=0\}=e^{-\lambda
\Delta_n}=0$ the result \eqref{fisher} follows.$\hfill\square$
\\

Fisher's information plays a central role in the Cram\'{e}r-Rao
inequality and more in general in the parametric inference. Let
$\mathbf{E}_\lambda^n(\cdot)$ be the expectation with respect to
the probability measure $\mathbf{P}_\lambda^n$. For any estimators
of the parameter $\lambda$, we have the next result.
\begin{theorem}
Let $T_n$ be an arbitrary estimator of $\lambda$ such that
$\mathbf{E}_\lambda^n |T_n|^2< \infty$ for any $\lambda>0$. Then
\begin{equation}
b(\lambda)=\mathbf{E}_\lambda^n T_n-\lambda \end{equation} is
differentiable respect to $\lambda$ in $\mathbf{L}_2$. Moreover,
the following Cram\'{e}r-Rao inequality holds
 \begin{equation}\label{cap.4:cr}
\mathbf{E}_\lambda^n(T_n-\lambda)^2\geq \frac{(1+\frac{d
b(\lambda)}{d\lambda})^2}{I_n(\lambda)}+b^2(\lambda).
\end{equation}
\end{theorem}
{\bf Proof.} We note that
\begin{equation}
\mathbf{E}_\lambda^nT_n=\mathbf{E}_{\lambda_0}^n\left\{T_n
\left(\frac{\lambda}{\lambda_0}\right)^n\exp\left(-(\lambda-\lambda_0)n\Delta_n
+\frac{\lambda-\lambda_0}{c}\sum_{i=1}^n\sqrt{u_{n,i}}\right)\right\}
\end{equation}
and show that $\mathbf{E}_\lambda^n T$ is differentiable and the
equality
\begin{equation}\label{cap.4:diff.cr}
\frac{d}{d \lambda}\mathbf{E}_\lambda^n
T_n=\mathbf{E}_\lambda^n\left\{T_n
\left(\frac{n}{\lambda}-n\Delta_n+\frac{1}{c}\sum_{i=1}^n\sqrt{u_{n,i}}\right)\right\}
\end{equation}
holds in $\mathbf{L}_2$.

For this purpose we interpret $\frac{d}{d
\lambda}\mathbf{E}_\lambda^n T_n$ as the right-hand side of
equation \eqref{cap.4:diff.cr}. It is not difficult to see that
\begin{eqnarray}\label{maj}
&&\Bigg|\mathbf{E}_{\lambda+h}^n T_n-\mathbf{E}_\lambda^n T_n - h
\frac{d}{d
\lambda}\mathbf{E}_\lambda^n T_n\Bigg|^2\\
&&=\Bigg|
\mathbf{E}_\lambda^n\left\{T_n\left[\frac{d\mathbf{P}_{\lambda+h}^n}{d\mathbf{P}_\lambda^n}-1-h\left(\frac{n}{\lambda}-n\Delta_n+\frac{1}{c}\sum_{i=1}^n\sqrt{u_{n,i}}\right)\right]\right\}
\Bigg|^2\notag\\
&&\leq \mathbf{E}_\lambda^n |T_n|^2
\mathbf{E}_\lambda^n\left\{e^{-h
n\Delta_n+\frac{h}{c}\sum_{i=1}^n\sqrt{u_{n,i}}+n\log(1+h/\lambda)}-1-h\left(\frac{n}{\lambda}-n\Delta_n+\frac{1}{c}\sum_{i=1}^n\sqrt{u_{n,i}}\right)
\right\}^2,\notag
\end{eqnarray}
where in the last step we have used the Cauchy-Schwarz inequality.

By inserting the equality
\begin{equation*}
e^{-h n
\Delta_n+\frac{h}{c}\sum_{i=1}^n\sqrt{u_{n,i}}+n\log(1+h/\lambda)}=1+h\left(\frac{n}{\lambda}-n\Delta_n+\frac{1}{c}\sum_{i=1}^n\sqrt{u_{n,i}}\right)
+o(h),
\end{equation*}
into \eqref{maj}, we can conclude that $b(\lambda)$ is
differentiable in $\mathbf{L}_2$-sense.

The validity of the inequality \eqref{cap.4:cr} follows by
standard arguments. \hfill $\square $\\

{\bf Remark 3.3} By taking into account an unbiased estimator
$T_n$ of the parameter $\lambda$, from \eqref{cap.4:cr} we get
that
\begin{equation}
\mathbf{E}_\lambda^n(T_n-\lambda)^2\geq \frac{\lambda^2}{n}.
\end{equation}

It's well-known that the Cram\'{e}r-Rao doesn't give a good
definition of asymptotic efficiency, because the limit variance
may not coincide with the variance of the limiting distribution.
Therefore to investigate the asymptotic properties of the
estimator $\widetilde{\lambda}_n$ as $n \rightarrow \infty$ and
$\Delta_n$ fixed, we reduce our problem to the study of the
normalized pseudo-likelihood ratio
\begin{eqnarray}\label{normlik}
Z_{n,\lambda}(z)&=&\frac{d\mathbf{P}_{\lambda+\varphi(n)z}^n}{d\mathbf{P}_{\lambda}^n}=\\
&=& \prod_{i=1}^n\exp\left(
\frac{\varphi(n)z}{c}\sqrt{u_{n,i}}-\varphi(n)z\Delta_n+\log\left(\frac{\lambda+\varphi(n)z}{\lambda}\right)\right)\notag\\
&=&\exp\left(
\frac{\varphi(n)z}{c}\sum_{i=1}^n\sqrt{u_{n,i}}-\varphi(n)nz\Delta_n+n\log\left(\frac{\lambda+\varphi(n)z}{\lambda}\right)\right),\notag
\end{eqnarray}
where $\varphi(n)=\varphi(n,\lambda)=(I_n(\lambda))^{-1/2}$. The
function \eqref{normlik} takes values in the following set
\[
U_{n,\lambda}=\left\{z:\lambda+\frac{z}{\sqrt{I_n(\lambda)}}\in
\Theta\right\}.
\]

It's well-known that $Z_{n,\lambda}$ (deriving from an i.i.d. observation scheme) admits the representation
\begin{eqnarray}
Z_{n,\lambda}(z)=\exp\left\{\frac{z}{\sqrt{I_n(\lambda)}}\sum_{i=1}^n\frac{\partial
\log
p(\lambda)}{\partial\lambda}-\frac{|z|^2}{2}+\phi_n(z,\lambda)\right\},
\end{eqnarray}
with $\frac{1}{\sqrt{I_n(\lambda)}}\sum_{i=1}^n\frac{\partial
p(\lambda)}{\partial\lambda}\overset{d}{\rightarrow} N(0,1)$ and
$\phi_n(z,\lambda)\rightarrow 0$ in probability as $n\rightarrow
\infty$; i.e. $\mathbf{P}_\lambda^n$ is locally asymptotically
normal (LAN).

 For the function $Z_{n,\lambda}$ we have the next useful
Lemma.
\begin{lemma}\label{cap.4:N3}
Let $K$ be a compact subset of $\Theta$. We have that:
\begin{itemize} \item[i)] for some constant $a=a(K), B=B(K)$
\begin{eqnarray}\label{con1}
\sup_{\lambda\in
K}\sup_{|z|<R,|v|<R}|z-v|^{-2}\mathbf{E}_\lambda^n\left|Z_{n,\lambda}^{1/2}(z)-Z_{n,\lambda}^{1/2}(v)\right|^2<B(1+R^a),
\end{eqnarray}
 with $z,v \in U_{n,\lambda};$\\
\item[ii)] for any $z\in U_{n,\lambda}$
\begin{eqnarray}\label{con2}
\sup_{\lambda\in K}\mathbf{E}_\lambda^n Z_{n,\lambda}^{1/2}(z)\leq
e^{-c|z|^2},
\end{eqnarray}
where $c>0$.
\end{itemize}
\end{lemma}
{\bf Proof.} $i)$ Following the proof of Lemma 1.1, section III in
Ibragimov-Has'minskii (1981) we get that
\begin{eqnarray}\label{lemn3}
&&\mathbf{E}_\lambda^n\left|Z_{n,\lambda}^{1/2}(z)-Z_{n,\lambda}^{1/2}(v)\right|^2\\
&&\leq\left|(I_n(\lambda))^{-1}\int_0^1I_n(\lambda+\varphi(n)(z+s(v-z))))ds\right||z-v|^2.\notag
\end{eqnarray}
 Since
 \begin{eqnarray*}
\frac{I(\lambda+z)}{
I(\lambda)}&=&\left(\frac{\lambda}{\lambda+z}\right)^2,
\end{eqnarray*}
for $z>0$ follows $\left(\frac{\lambda}{\lambda+z}\right)^2\leq
1$. If $z<0$ we can see that
\begin{eqnarray*}
\left(1+\frac{z}{\lambda}\right)^{-2}=1+(-2)\frac{z}{\lambda}+o\left(\frac{z}{\lambda}\right)<3+o\left(\frac{z}{\lambda}\right).
\end{eqnarray*}

Therefore set $B=3+o\left(\frac{z}{\lambda}\right)$ the inequality
\begin{equation}\label{fishratio}
\sup_{\lambda\in\Theta}\sup_{|z|<R,\lambda+z\in
\Theta}\left|\frac{I_n(\lambda+z)}{I_n(\lambda)}\right|\leq
B(1+R^a),
\end{equation}
holds.

In view of the relationships \eqref{fishratio} and \eqref{lemn3}
the proof of the inequality \eqref{con1} is concluded.

$ii)$ The function $\partial\sqrt{p(\lambda)}/\partial\lambda$ is
differentiable in $\mathbf{L}_2$, then
\begin{eqnarray*}
\iint_{S_{c\Delta_n}^2}|\sqrt{p(\lambda+h)}-\sqrt{p(\lambda)}|^2dxdy&=&\iint_{S_{c\Delta_n}^2}h^2
\left(\frac{\partial}{\partial\lambda}\sqrt{p(\lambda)}
\right)^2dxdy+o(|h|^2)\\
&=&\frac{h^2}{4}I(\lambda)+o(|h|^2).
\end{eqnarray*}

By taking into account that
$$0<\inf_{\lambda \in \Theta}I(\lambda)< \sup_{\lambda \in \Theta}
I(\lambda)< \infty,$$ we have immediately that

\begin{equation}\label{pos}
\iint_{S_{c\Delta_n}^2}|\sqrt{p(\lambda+h)}-\sqrt{p(\lambda)}|^2dxdy>0.
\end{equation}

From \eqref{pos}, we derive the inequality
\begin{eqnarray*}
\inf_{\lambda \in K}\inf_{\{ h:\lambda+h\in
\Theta\}}\iint_{S_{c\Delta_n}^2}|\sqrt{p(\lambda+h)}-\sqrt{p(\lambda)}|^2
dxdy\geq \frac{a|h|^2}{1+|h|^2},\quad a>0,
\end{eqnarray*}
and Lemma 5.3, Chapter I, in Ibragimov and Has'minskii (1981)
permits us to obtain the condition
\eqref{con2}.$\hfill\square$\\

Finally, we are able to present the main result of this section.

\begin{theorem}\label{main}
Let $K$ be a compact subset of $\Theta$. Then the estimator
$\widetilde{\lambda}_n$, defined in \eqref{cap.4:maxlikmod},
uniformly in $\lambda\in K$:
\begin{itemize} \item is consistent; \item converges in
distribution as follows
\begin{equation}
\sqrt{I_n(\lambda)}(\widetilde{\lambda}_n
-\lambda)\overset{d}{\rightarrow} N(0,1);
\end{equation}
\item has moments such that
\begin{equation}
\lim_{n \rightarrow
\infty}\mathbf{E}_\lambda^{n}|\sqrt{I_n(\lambda)}(\widetilde{\lambda}_n
-\lambda)|^\gamma=\mathbf{E}|\xi|^\gamma,
\end{equation}
 where $\gamma>0$ and $\xi \sim
N(0,1)$.
\end{itemize}
\end{theorem}
{\bf Proof.} In accordance with the Theorem 1.1 Chapter III, in
Ibragimov and Has'minskii (1981), we prove that the four
conditions are satisfied.

The probability measure $\mathbf{P}_\lambda^n$ is uniformly local
asymptotic normal, while it's easy to see that
$$\lim_{n\rightarrow \infty}\sup_{\lambda\in K}\varphi^2(n,\lambda)=0.$$

The validity of Lemma \ref{cap.4:N3} concludes the proof.$\hfill\square$\\

The Theorem \ref{main} yields the H\'{a}jek-Le Cam asymptotic
efficiency of the estimator $\widetilde{\lambda}_n$ with respect
to a quadratic loss function. In fact, we have that
\begin{equation}
\lim_{\delta\rightarrow0}\varliminf_{n\rightarrow
\infty}\sup_{|\lambda-\lambda_0|<\delta}
\mathbf{E}_\lambda^{n}\left|\sqrt{I_n(\lambda)}(\widetilde{\lambda}_n-\lambda)\right|^2=1.
\end{equation}

\section{Monte Carlo analysis}
We analyze the empirical performance of the pseudo-maximum
likelihood estimator $\hat{\lambda}_n$ by means of a Monte Carlo
analysis with $n<\infty$ fixed. We simulate 10000 sample paths of
the planar random flights in the interval $[0,T]$, with $T=500$,
for different values of $\lambda$ and $c=1$. For any trajectories
we have sampled $n=200,300,500,1000$ values subsequently used to
estimate the unknown parameter $\lambda$.

The results have been reported in the Table \ref{cap.3:tab1.bis}.
Furthermore in the Table \ref{cap.3:tab1.bis} there is a column
$\sqrt{\text{MSE}( \lambda)}$ derived as follows
\begin{equation} \sqrt{\text{MSE}(
 \lambda)}=\sqrt{\frac{1}{N}\sum_{i=1}^N(\hat{\lambda}_n-\lambda)^2},
\end{equation}
where $N=10000$ is the number of simulations.

 It emerges, as expected, that the mean square error
  tends to zero when the sample size increases. Furthermore, it is clear that the true value
   of the parameter $\lambda$ and the mean square error are
   correlated. In fact, for fixed $n$, as the more $\lambda$ increases the more Poisson events remain
   hidden to the observer. The bias assumes small values for all the cases considered and
   is constantly equal to 0.002 for
$\lambda=0.1,\,0.25,\,0.5,\,0.75$.
\begin{table}[ht]
\begin{center}
\begin{tabular}{r|rrrr|r}
\hline
 $\lambda$ & Bias & $\sqrt{{\rm MSE}(\lambda)}$ & $\min\hat\lambda_n$ & $\max\hat\lambda_n$ & $n$ \\
\hline
 0.10 & 0.002 & 0.015 & 0.05 & 0.15 &  200 \\
 & 0.002 &0.015 & 0.05 & 0.16 & 300 \\
 & 0.002 & 0.015 & 0.05 & 0.15 & 500 \\
 & 0.002 & 0.014 & 0.06 & 0.15 & 1000 \\
 0.25 & 0.002 & 0.026 & 0.17 & 0.37 &  200 \\
 & 0.002 & 0.025 & 0.17 & 0.35 & 300 \\
 & 0.002 & 0.024 & 0.17 & 0.35 & 500 \\
 & 0.002 & 0.023 & 0.17 & 0.34 & 1000 \\
 0.50 & 0.001 & 0.042 & 0.37 & 0.67 &  200 \\
 & 0.002 & 0.038 & 0.36 & 0.65 & 300 \\
 & 0.002 & 0.035 & 0.36 & 0.65 & 500 \\
 & 0.002 & 0.033 & 0.36 & 0.63 & 1000 \\
 0.75 & $-$0.000 & 0.057 & 0.56 & 1.05 &  200 \\
 & 0.001 & 0.051 & 0.56 & 0.98 & 300 \\
 & 0.002 & 0.046 & 0.60 & 0.99 & 500 \\
 & 0.002 & 0.042 & 0.62 & 0.93 & 1000 \\
 1.00 & $-$0.004 & 0.073 & 0.76 & 1.28 &  200 \\
& $-$0.001 & 0.064 & 0.76 & 1.29 & 300 \\
 & 0.001 & 0.056 & 0.81 & 1.26 & 500 \\
 & 0.002 & 0.050 & 0.82 & 1.18 & 1000 \\
 1.50 & $-$0.013 & 0.106 & 1.15 & 1.98 &  200 \\
 & $-$0.003 & 0.090 & 1.19 & 1.92 & 300 \\
 & 0.001 & 0.076 & 1.22 & 1.78 & 500 \\
 & 0.001 & 0.066 & 1.26 & 1.78 & 1000 \\
 2.00 & $-$0.031 & 0.141 & 1.49 & 2.53 &  200 \\
 & $-$0.010 & 0.117 & 1.57 & 2.61 & 300 \\
 & 0.000 & 0.097 & 1.67 & 2.41 & 500 \\
 & 0.001 & 0.080 & 1.69 & 2.29 & 1000 \\
\hline\hline
\end{tabular}
\end{center}
\caption{Empirical performance of the estimator $\hat\lambda_n$
defined in \eqref{cap.4:maxlikmod} for different values of the
parameter $\lambda$ and different sample sizes. The velocity $c$
assumes value 1. The time horizon $T$ is equal to 500. The results
have been obtained on 10000 Monte Carlo sample paths of the planar
random flights. } \label{cap.3:tab1.bis}
\end{table}

\end{document}